\documentclass[draft]{article}

\usepackage{amsmath,amsthm}

\def\Int{\operatorname{Int}}
\def\diam{\operatorname{diam}}

\def\wt{\operatorname{wt}}
\def\lg{\operatorname{lg}}
\def\bord{\partial}

\let\bydef\emph

\def\var{ma\-ni\-fold}
\def\svar{sub\-ma\-ni\-fold}

\newtheorem{theo}{Theorem}[section]

\newtheorem{prop}[theo]{Proposition}
\newtheorem{lem}[theo]{Lemma}

\theoremstyle{definition}

\newtheorem*{claim}{Claim}

\newenvironment{rems}{\bigskip\emph{Remarks.}\begin{itemize}}{\end{itemize}}
\newenvironment{concrems}{\bigskip\emph{Concluding remarks.}\begin{itemize}}{\end{itemize}}

\def\Rr{\mathbf{R}}

\def\Nn{\mathbf{N}}

\def\Ss{\mathbf{S}}

\def\calF{\mathcal{F}}
\def\calS{\mathcal{S}}
\def\calT{\mathcal{T}}
\def\calD{\mathcal{D}}

\begin{document}
\title{A spherical decomposition for Riemannian open $3$-manifolds}
\author{Sylvain Maillot}
\date{Final version, May 2006}

\maketitle

\begin{abstract}
We show that open $3$-manifolds that have a locally finite
decomposition along $2$-spheres are characterized
by the existence of a Riemannian
metric with respect to which the second homotopy group of the manifold
is generated by small elements.
\end{abstract}

\section{Introduction}
Open $3$-manifolds have been studied from the topological viewpoint
for a long time. A wealth of examples of `exotic' open $3$-manifolds
(see e.g.~\cite{whitehead:unity,scott:noncompact,
st:exotic}) show that many classical results in $3$-manifold topology
do not hold when the manifold is not assumed to be compact.
By contrast, much work has been done towards proving that such exotic
behaviors do not occur in the realm of \emph{hyperbolic} open $3$-manifolds,
leading to the recent proof of the Marden Conjecture by Calegari-Gabai,
and independently Agol
(see~\cite{agol:tame,cg:tame} and the references therein.)

In this paper, we study open $3$-manifolds
endowed with Riemannian metrics controlling the topology on the large
scale, in the spirit of Gromov~\cite{glp:struc},
without imposing a sign on the curvature.
Our goal is to find an extension of Kneser's Prime Decomposition
Theorem~\cite{kneser:sphere} to this context. An example of
P.~Scott~\cite{scott:noncompact} shows that an open $3$-manifold need
not have a prime decomposition, even allowing infinitely many factors.
The reason is that a maximal collection of disjoint essential spheres
may fail to be \emph{locally finite.}

To state our main result, we need some terminology.  A $3$-manifold
$M$ is called \bydef{weakly irreducible} if the manifold $\hat M$
obtained from $M$ by capping off $3$-balls to all $2$-sphere boundary
components has the following property: every embedded $2$-sphere in
$\hat M$ bounds a compact contractible submanifold.  A Riemannian
\var\ has \bydef{bounded geometry} if the absolute values of its
sectional curvatures are uniformly bounded and its injectivity radius
is positive.

\begin{theo}[Main Theorem]\label{main}
Let $M$ be an orientable $3$-manifold without boundary.  Suppose that
there exist a complete Riemannian metric $g$ of bounded geometry on $M$ and a
constant $C\ge0$ such that $\pi_2M$ is generated as a $\pi_1M$-module by
homotopy classes of smooth maps from the $2$-sphere into $M$ whose
area, measured with respect to $g$, is at most $C$. Then there is a
locally finite collection $\calS$ of pairwise disjoint embedded $2$-spheres in
$M$ that splits $M$ into weakly irreducible \svar s.
\end{theo}

\begin{rems}
\item The definition of weak irreducibility has been phrased so as to
make the statement of the Main Theorem independent
of the Poincar\'e Conjecture. If this
conjecture is true, then the conclusion of the Main Theorem is a
natural generalization of the existence of prime decompositions for
compact manifolds. Note that we allow nonseparating spheres in
$\calS$; those can be replaced by $S^2\times S^1$ factors if one
wishes to use the language of connected sum decompositions.
\item If a $3$-manifold $M$ satisfies the conclusion of the above
theorem, then it is easy to construct a complete Riemannian metric of bounded
geometry on $M$ such that all $2$-spheres in $\calS$ have area $1$. It
follows from Proposition~\ref{top equiv} that $\calS$ generates
$\pi_2M$ as a $\pi_1M$-module. Thus the sufficient condition of the
Main Theorem is also necessary.
\item In the statement of the hypothesis, a standard abuse of language
has been made: strictly speaking, classes in $\pi_2M$ are represented,
not by maps $S^2\to M$, but by based maps $(S^2,*)\to (M,*)$. Once
basepoints have been chosen, any map can be homotoped to a based map.
The ambiguity in doing this is measured by the action of $\pi_1M$ on
$\pi_2M$. In the sequel, all statements will be invariant under this action;
hence reference to basepoints may safely be dropped.
\item There is no direct connection between our work and the issue of
\emph{topological tameness,} i.e.~finding conditions for an open \var\
with finitely generated fundamental group to be homeomorphic to the
interior of some compact manifold.  The $3$-\var s we are interested in
typically have huge second homotopy group, and their fundamental
groups may or may not be finitely generated.
\end{rems}

Before closing this introduction, I would like to thank Steven Boyer,
Peter Shalen and Ian Agol for discussions, and acknowledge financial
support from the CRM and the CIRGET at Montreal, where a substantial
part of this work was done.

\section{Topological preliminaries} 
In this section, we prove a purely topological criterion for a collection
of spheres to split $M$ into weakly irreducible \svar s. For convenience,
we introduce some terminology and notation which will be used throughout
the paper.

A \bydef{system} of surfaces in a $3$-manifold $M$ is a locally finite
collection $\calF$ of pairwise disjoint surfaces embedded in $M$. Those
surfaces are called the \bydef{components} of $\calF$. If all
components are spheres, we say that $\calF$ is \bydef{spherical}.

If $\calF$ is a system of surfaces, we shall denote by
$M\backslash\calF$ the manifold obtained from $M$ by removing a
disjoint union of open product neighborhoods of the components of
$\calF$. The operation of removing such neighborhoods is called
\bydef{splitting} $M$ \bydef{along} $\calF$.

\begin{prop}\label{top equiv}
Let $M$ be an orientable $3$-manifold with empty boundary
and $\calS$ a spherical system
in $M$. Then $\calS$ generates $\pi_2M$ as a $\pi_1M$-module
if and only if each component of $M\backslash\calS$ is weakly irreducible.
\end{prop}

\begin{proof}
First we prove the `if' part. Assume that each component of
$M\backslash\calS$ is weakly irreducible. Seeking a contradiction,
suppose that the $\pi_1M$-submodule $A$ of $\pi_2M$ generated by
$\calS$ is proper. Then by the Sphere Theorem, there is an embedded
$2$-sphere $S$ in $M$ whose homotopy class does not belong to $A$.

Assume that $S$ is in general position with respect to $\calS$
and intersects it in the least possible number of components. If
$\bigcup\calS\cap S$ is nonempty, then we can find a disc
$D\subset\bigcup\calS$ bounded by a curve of $\bigcup\calS\cap S$ and
whose interior does not intersect $S$. Surgering $S$ along $D$, we get
two $2$-spheres, one of which at least does not belong to $A$, and has
fewer intersection components with $\calS$ than $S$, a contradiction.

Hence $S$ lies in a component $X$ of $M\backslash\calS$. Let $\hat X$ be
the $3$-manifold obtained from $X$ by capping off $3$-balls. By
hypothesis, $S$ bounds a contractible compact submanifold of $\hat
X$. Therefore, $M$ contains a compact simply-connected submanifold $Y$ such
that $\bord Y$ is the union of $S$ and finitely many components of
$\calS$. This means that $S$ belongs to the submodule generated by the
other components of $\bord Y$, hence to $A$, contradicting our hypothesis.
This proves the `if' part. 

Let us turn to the `only if' part. Let $\cal S$ be a spherical system
generating $\pi_2M$ as a $\pi_1M$-module. Let $X$ be the $3$-complex
obtained from $M$ by gluing a $3$-ball onto each sphere in $\cal S$.
Then the universal cover $\tilde X$ of $X$ is $2$-connected.
By the Hurewicz Isomorphism Theorem, $H_2(\tilde X)=0$.

Let $Y$ be a component of $M \backslash \cal S$ and let $N$ be obtained
from $Y$ by capping off a $3$-ball to each boundary component of $Y$.
Then $\tilde N$ embeds into $\tilde X$ in an obvious way. Using the
Mayer-Vietoris exact sequence, it follows that $H_2(\tilde N)=0$.
By the Hurewicz Isomorphism Theorem again, $\pi_2 \tilde N=\pi_2 N=0$.
Hence any embedded $2$-sphere in $N$ bounds a contractible manifold in $N$.
\end{proof}

\section{Normal surfaces in Riemannian 3-manifolds}

In all of this subsection, let $(M,\calT)$ be an orientable $3$-manifold
with a fixed triangulation. 

The \bydef{size} of a subset $A\subset M$ is the minimal number of
$3$-simplices of $\calT$ needed to cover $A$. Then we define a
quasimetric (see~\cite{sm:seifert}) $d_\calT$ on $M$ as follows:
given two points $x,y\in M$, we let $d_\calT(x,y)$ be the minimal size
of a path connecting $x$ to $y$ minus one. (Quasi)metric balls,
neighborhoods, diameter etc.~can be defined in the usual way.

Recall from~\cite{sm:seifert} the definition of a \bydef{regular
Jaco-Rubinstein metric} on $(M,\calT)$: it is a Riemannian metric on
$\calT^{(2)}-\calT^{(0)}$ such that each $2$-simplex is sent
isometrically by barycentric coordinates to a fixed ideal triangle in
the hy\-per\-bo\-lic plane. The crucial property for
applications to noncompact manifolds is that for every number $n$,
there are finitely many subcomplexes of size $n$ up to isometry.

Let $F$ be a compact, orientable
surface and $f:F\to M$ be a proper map in general position with
respect to $\calT$. The \bydef{weight} $\wt(f)$ of $f$ is the number of
points of $f(F)\cap \calT^{(1)}$ counted with multiplicities.
Its \bydef{length} $\lg(f)$ is the total length of all the arcs in the
boundaries of the disks in which $f(F)$ intersects the 3-simplices of
$\mathcal T$. The \bydef{PL area} of $f$ is the pair 
$|f|=(\wt(f),\lg(f))\in \Nn\times \Rr_+$. We are interested in surfaces
having least PL area among surfaces in a particular class with 
respect to the lexicographic order.

We say that $f$ is
\bydef{normal} if $f(F)$ misses $\calT^{(0)}$ and meets
transversely each 3-simplex $\sigma$ of $\calT$ in a finite collection of
disks that intersect each edge of $\sigma$ in at most one point.

The following lemma from~\cite{sm:seifert} provides a useful
inequality between the weight of a normal surface and the diameter of
its image with respect to the quasimetric $d_\calT$. The statement
given here combines Lemmas~2.2 and~A.1 of~\cite{sm:seifert}.

\begin{lem}\label{myineq}
Let $F$ be a compact surface and $f:F\to M$ be a normal map. Then
$\diam(f(F))\le \wt(f)^2$.
\end{lem}

We will face repeatedly the
following problem:

(*) Let $A$ be a proper submodule of
$\pi_2M$ (viewed as a $\pi_1M$-module). Find a map $f\!:\Ss^2\to M$
whose homotopy class does not belong to $A$, and
with least PL area among such maps.

\begin{lem}\label{minimizers}
Let $A$ be a proper submodule of $\pi_2M$.
\begin{enumerate}
\item Problem (*) has a solution.
\item If $f_0\!:\Ss^2\to M$ is such a solution, then $f_0$ is normal;
furthermore, $f_0$ is either an embedding or a double cover of an
embedded projective plane.
\item If $B$ is another proper submodule of $\pi_2M$, $f$ (resp.~$g$)
a solution of (*) for $A$ (resp.~$B$), then the images of $f$ and $g$
are disjoint or equal.
\end{enumerate}
\end{lem}

\begin{proof}
Proof of (i): let ${\cal S}_A$ be the class of all maps $f\!:\Ss^2\to
M$ whose homotopy class does not belong to $A$. Since the weight is a
nonnegative integer, there is a least weight map $f\in{\cal S}_A$.
Applying the Kneser-Haken normalization process, we find a least
weight \emph{normal} map $f'\in{\cal S}_A$ with $|f'| \le |f|$. It is
now sufficient to find a map of least length in the class ${\cal
S}^0_A$ of least weight normal maps in ${\cal S}_A$. Such a minimizer
exists, because by lemma~\ref{myineq} there are only finitely many
combinatorial types for members of ${\cal S}^0_A$, so by regularity of
the Jaco-Rubinstein metric, the set of lengths of members of ${\cal
S}^0_A$ is finite.

Proof of (ii): Let $f_0$ be a minimizer. Since the normalization
process strictly decreases PL area, $f_0$ must be normal. In
particular it has least PL area in its normal homotopy class.

To prove embeddedness, we will use a trick described in the Appendix
of~\cite{jr:min}, which consists in perturbing the Jaco-Rubinstein
metric in a neighborhood of the image of $f_0$ to achieve general
position for $f_0$. The perturbed metric will not be regular, but will still
have the property that for every $n$,
there are finitely many subcomplexes of size $n$ up to isometry, so
noncompactness of $M$ will not be a problem.
Hence Corollary~3 and Remark~7 of~\cite{jr:min} apply. As a consequence,
it suffices to show that $f_0$ can be perturbed to a map
which is an embedding or a double cover.

We perform the perturbation so that the resulting map, which we still
call $f_0$, is in general position.  Then we show how to adapt
arguments of Meeks-Yau~\cite[Theorem~6]{my:embedding}. The first
part of the proof consists in constructing a tower of coverings of
some regular neighborhood $N$ of the image of $f_0$. This part is
topological and goes through without changes.

The crucial part where
the geometry is used is in Assertion~1 on page~471 of~\cite{my:embedding},
which states that the lift
$f_k:\Ss^2 \to N_k$ of $f_0$ to the top of the tower is an embedding.
To prove this, one compares the area of $f_k$ with the areas of
various spheres which are components of $\bord N_k$. One argues that
if the image of $f_k$ were not embedded, these spheres would have
folding lines allowing their area to be reduced by rounding off; 
however one of them must be topologically essential when projected
down the tower, giving a contradiction.

To adapt the argument to the PL setting, one only needs a coherent
way of measuring PL areas in $N_k$. To do this, observe that
$\calT$ induces a cell decomposition of $N$ (which may not be
a triangulation). This cell decomposition can be lifted to $N_k$
and used to measure PL area. Then the Meeks-Yau argument can be
adapted using the work of Jaco-Rubinstein. Once one knows
that $f_k$ is an embedding, the rest of the proof goes through
without changes.

Proof of (iii): using the same perturbing trick (cf.~Corollary~4 and
Remark~7 of~\cite{jr:min}), we may assume that $f,g$
are embeddings and their images $S_1,S_2$ intersect transversely. Let
$D$ be a disk embedded in $S_1$ or $S_2$, whose boundary is an
intersection curve. Assume that $D$ has least PL area among such
disks. Then the interior of $D$ contains no curve of $S_1\cap
S_2$. Without loss of generality we assume that $D\subset D_1$. Then
$S_2\cap D$ splits $S_2$ into two disks $D',D''$.

Let $S'$ (resp.~$S''$) be obtained from $D\cup D'$ (resp.~$D\cup D''$) by
rounding the corners so that $|S'| < |D|+|D'|$ and $|S''| <
|D|+|D''|$. It follows from the minimality hypothesis of $D$ that $S'$
and $S''$ both have PL area strictly less than $S_2$. Now at least one
of them is not in $B$, contradicting the hypothesis on $S_2$.
\end{proof}

The connection between Riemannian geometry and PL topology is provided
by the following lemma.

\begin{lem}\label{push}
Let $g$ be a Riemannian metric on $M$ such that every $3$-simplex of
$\calT$ is bi-Lipschitz homeomorphic to the standard $3$-simplex with
uniform Lipschitz constants.
Then there is a constant $\delta>0$ such that for every smooth map
$f\!:\Ss^2\to M$ there are smooth maps
$f_1,\ldots,f_n\!:\Ss^2\to M$ such that $f$ is contained in the
$\pi_1M$-submodule of $\pi_2M$ generated by $f_1,\ldots,f_n$ and for
each $i$, the weight of $f_i$ is bounded above by $\delta$ times the
$g$-area of $f$.
\end{lem}

\begin{proof}
Let $f\!:\Ss^2\to M$ be a smooth map. We are going to perform a number of
modifications on $f$, checking at each stage that the $g$-area does
not increase by more than a multiplicative factor. To keep notation
simple, we will still denote the resulting map by $f$. In the final
step, the maps $f_1,\ldots,f_n$ will appear.

\paragraph{Step 1} Modify $f$ so that $f(\Ss^2)\subset\calT^{(2)}$. 

The needed argument is essentially contained in~\cite{bt:dehn}, but
we reproduce it here for completeness.

We know that $3$-simplices of our triangulation are uniformly bi-Lipschitz
equivalent to the standard simplex $\sigma_0$. For each $3$-simplex $\sigma$ we
choose once and for all a bi-Lipschitz parametrisation
$\phi_\sigma:\sigma\to\sigma_0$. If $p$ is a point of $\Int\sigma_0$,
let $\psi_p\!:\sigma_0-\{p\}\to\bord\sigma_0$ denote the radial projection centered
at $p$. If $\sigma$ is a $3$-simplex of $\calT$,
the \bydef{radial projection centered at a point} $p\in\Int\sigma$ is
defined as $\phi_\sigma \circ \psi_{\phi_\sigma(p)} \circ \phi_\sigma^{-1}$.
We shall prove the following:

\begin{claim}
There is a constant $\lambda>0$ depending on $\calT$ and $g$, but not
on $f$, so that one can choose for each $3$-simplex $\sigma$ a point
$p_\sigma\in\Int \sigma-(\sigma\cap f(\Ss^2))$ depending on $f$
so that composing $f$ with the
radial projections $\sigma\to\bord\sigma$ centered at $p_\sigma$
can increase area by at most a multiplicative factor $\lambda$.
\end{claim}

To prove the claim, we consider a $3$-simplex $\sigma$. For simplicity
we assume that $\sigma$ is isometric to the standard (Euclidean)
$3$-simplex. Since the simplices of our triangulation are uniformly bi-Lipschitz
equivalent to the standard simplex, this induces no loss of generality
(only the values of the constants involved are changed). Let $0$ be
the barycenter of $\sigma$, $r>0$ be a constant such that
$B(0,3r)\subset\Int\sigma$. Set $B:=B(0,r)$ and $Q:=f(\Ss^2)$.
For every $u\in B$,
let $B_u$ denote the ball around $y$ of radius $2r$. By hypothesis, we
have $B\subset B_u\subset\sigma$. Let $\pi_u\!:B_u-\{u\}\to\bord B_u$
denote the radial projection centered at $u$. By convention, we extend
$\pi_u$ as the identity on $\sigma-B_u$.

Let $|X|_i$ denote the $i$-dimensional Hausdorff measure of a
subset $X\subset \sigma$.

Radial projections have the property that away from a ball of given
radius $\rho$ around the center of the projection, the increase in
area is bounded above by a multiplicative factor that depends only on
$\rho$.  Hence if we were able to find a point $u$ whose distance to
$Q$ is bounded below independently of $f$, the claim would follow. Of
course this need not be true in general, but we shall find $u$ such
that the restriction of $\pi_u$ to $Q$ does not increase area too
much. The bound of the area dilatation of ${\psi_{u}}_{|Q}$ will
follow.

If $|Q|_2=0$, there is nothing to prove, so assume $|Q|_2\neq0$.
Arguing by contradiction, we consider for every $\nu>0$ the `bad' set $A_\nu$ of
points $u\in B$ such that $u\not\in Q$ and $|\pi_u(Q)|_2>\nu|Q|_2$. Our
next goal is to derive some volume estimates for $A_\nu$. Since $\pi_u$
is the identity away from $B_u$, we have:
\begin{align*}
|\pi_u(Q)|_2&\le |Q\cap (\sigma-B_u)|_2 + |\pi_u(Q\cap B_u)|_2\\
 &\le |Q|_2 + \int_{Q\cap B_u}\,\, {\left(\frac{2r}{\|x-u\|} \right)}^2\,dx.
\end{align*}

Hence
\begin{align*}
|A_\nu|_3 &= \int_{A_\nu}\,du = (\nu|Q|_2)^{-1} \cdot
 \int_{A_\nu}\,\,\nu|Q|_2\,du\\
&\le (\nu|Q|_2)^{-1} \cdot \int_{A_\nu} |\pi_u(Q)|_2\,du
 \le (\nu|Q|_2)^{-1} \cdot \int_{B} |\pi_u(Q)|_2\,du\\
& \le (\nu|Q|_2)^{-1} \cdot \left({\int_B\,\,|Q|_2\,du} + {\int_B\,\,du\,\int_{Q\cap B_u}\,\,
\frac{4r^2}{\|x-u\|^2}\,dx}\right)\\
& \le (\nu|Q|_2)^{-1} \cdot \left({|B|_3 \cdot |Q|_2} + {\int_B\,\,du\,\int_{Q\cap B_u}\,\,
\frac{4r^2}{\|x-u\|^2}\,dx}\right).
\end{align*}

Since $Q\cap B_u$ is a compact set and $\|x-u\|$ does not vanish on
this set, we can apply Fubini's theorem. After making the change of
variables $(y,z)=(x,x-u)$, we get:

\begin{align*}
\int_B\,\,du\,\int_{Q\cap B_u}\,\,\frac{4r^2}{\|x-u\|^2}\,dx
&=\int_Q\,\,dy\,\int_{B(x,r)\cap
B(0,2r)}\,\,\frac{4r^2}{\|z\|^2}\,dz\\
&\le |Q|_2 \cdot \int_{B(0,2r)}\,\,\frac{4r^2}{\|z\|^2}\,dz.
\end{align*}

Passing to polar coordinates, we see that the last integral is
bounded above by some constant $K$. We deduce

\begin{align*}
|A_\nu|_3 &\le(\nu|Q|_2)^{-1} \cdot \left(|B|_3 \cdot |Q|_2 +
K\cdot |Q|_2\right)\\
&\le \nu^{-1} (|B|_3 +K).
\end{align*}

Hence if we take $\nu_0$ large enough, we can make the volume of the bad
set $A_{\nu_0}$ as close to zero as we want. Since $f$ is smooth, $Q\cap
B$ has zero Lebesgue measure, so setting $\nu_0:=2(|B|_3 +K) \cdot
|B|^{-1}$ is sufficient to ensure that
$B-Q-A_{\nu_0}$ has nonzero  Lebesgue measure, hence is nonempty.
This proves the existence of a point $u\in B$ such that
$|\pi_u(Q)|_2\le \nu_0 \cdot |Q|_2$. As remarked before, this shows that
$|\psi_u(Q)|_2\le \lambda \cdot |Q|_2$ for some constant $\lambda$
depending only on $r$ and $\nu_0$, hence not on $Q$. This proves the
claim and completes our first step.

\paragraph{Step 2} Modify $f$ so that the image still lies in
$\calT^{(2)}$, and there is a triangulation $\calD$ of $\Ss^2$
that makes $f$ simplicial and a constant $\epsilon>0$ independent of
$f$ such that the number of $2$-simplices of
$\calD$ that are mapped homeomorphically to $2$-simplices of $\calT$
is bounded above by $\epsilon$ times the $g$-area of
$f$.

This is a standard argument, but some care is needed to get the
required upper bound. Note that there is no bound on the total number
of $2$-simplices in $\calD$, since our original map $f$ could have
$g$-area arbitrarily close to zero and yet be sent to a very long path
in $\calT^{(1)}$ by Step~1.

Let $\sigma$ be a $2$-simplex of $\calT$. The preimage by $f$ of
$\Int\sigma$ has finitely many components, which are open subsets of
$\Ss^2$. Let $X$ be one of them. If $f_{|X}$ is not onto, then we can
use radial projection from a point of $\Int\sigma$ to push it off
$\Int\sigma$. After finitely many such operations, we can assume that
$f_{|X}$ is onto for each component $X$ of $f^{-1}(\Int\sigma)$.

Then the degree $n_{\sigma,X}$
of the restriction $f:X\to\Int\sigma$ is bounded above by a
constant times the $g$-area of $f$. Let $q$ be a regular value of
this map and $p_1,\ldots,p_n$ be its preimage. Choose an open disk $V$
containing $q$ such that $f^{-1}(V)$ is a union of pairwise disjoint
open disks $U_1,\ldots,U_n$ with $p_i\in U_i$ for each $i$, and each
restriction $f\!:U_i\to V$ is a diffeomorphism.

Our next goal is to modify $f$ so that all $U_i$'s are mapped to $V$
with the same orientation. Assume that, say, $U_1$ and $U_2$ are
mapped with different orientations. Let $\xi$ be an arc in
$X-\bigcup_i U_i$ connecting $U_1$ to $U_2$. Since $f\!:X-\bigcup_i
U_i \to \sigma-V$ induces an epimorphism on fundamental groups, we can
choose $\xi$ so that $f\circ\xi$ is a null-homotopic loop in
$\sigma-V$.

Then $f$ can be homotoped so that $f\circ\xi$ is contracted to a
point. The number of components of $f^{-1}(V)$ decreases in the
process, and the area of $f$ does not increase. Hence after finitely
many of these modifications, $U_1,\ldots,U_n$ are all mapped to $V$
with the same orientation. In particular, $n=n_{\sigma,X}$.

Let $\tau$ be a triangle embedded in $V$. By composing $f$ with the
expansion of $\tau$ into $\sigma$, we modify it so that $f$ maps $n_{\sigma,X}$
disks $D_1(\sigma,X),\ldots,D_{n_{\sigma,X}}(\sigma,X)$
homeomorphically onto $\Int\sigma$ and the rest of $X$ to
$\calT^{(1)}$. Having done this for each $2$-simplex $\sigma$ and each
component $X$ of $f^{-1}(\Int\sigma)$, we choose a triangulation
$\calD$ of $\Ss^2$ such that the closures of the $D_i(\sigma,X)$ are
$2$-simplices. Then after a simplicial approximation on the part that
is mapped to $\calT^{(1)}$, $f$ is
simplicial with respect to $\calD$ and $\calT$, and the number of
$2$-simplices of $\calD$ that are mapped homeomorphically to
$2$-simplices of $\calT$ is bounded above by a multiplicative constant
times the $g$-area of $f$.

\paragraph{Step 3} The end.

We make an ordered list $\sigma_1,\ldots,\sigma_m$ of all
$2$-simplices of $\calD$ that are mapped into $\calT^{(1)}$ and
collapse them in that order. In the end, we get a finite cell complex with
the homotopy type of a bouquet of $2$-spheres, and the $2$-cells are
in bijection with the $2$-simplices of $\calD$ that are mapped
homeomorphically onto their images. Hence we can find combinatorial
maps $f_1,\ldots,f_n\!:\Ss^2\to \calT^{(2)}$ that generate $f$
homotopically and whose combinatorial areas are bounded by a constant
times the $g$-area of $f$.

Finally, for each $i$ we modify $f_i$ in
the following way: first push each $2$-simplex off the $2$-skeleton
into a neighboring $3$-simplex. This can be done consistently, and
hence realized by a homotopy on $f_i$, because of orientability. Then
approximate $f_i$ by a map in general position without introducing new
intersection points with the $1$-skeleton. Then we have the required
upper bound on the weight of $f_i$ and the proof of Lemma~\ref{push}
is complete.
\end{proof}

\section{Proof of the Main Theorem}
From now on, $M$ is an orientable $3$-manifold without boundary,
$g$ a complete Riemannian metric of bounded geometry on $M$, and
$C$ a constant
such that $\pi_2M$ is generated as a $\pi_1M$-module by homotopy classes of
spheres of area at most $C$. 

\begin{lem}\label{diam bound}
There exists a triangulation $\calT$ and a constant $C_1$ such that
for every proper submodule $A$ of $\pi_2M$, there is a solution to (*)
with diameter at most $C_1$.
\end{lem}

\begin{proof}
By~\cite{attie:bg}, the bounded geometry on $g$ implies that $M$ has
a triangulation $\calT$ such that there exists a constant $C_2>0$ such
that every $3$-simplex of $\calT$ is $C_2$-bi-Lipschitz homeomorphic to the
standard $3$-simplex. The fact that $\pi_2M$ is generated as a
$\pi_1M$-module by homotopy classes of spheres of $g$-area at most $C$
together with Lemma~\ref{push} implies that $\pi_2M$ is generated as a
$\pi_1M$-module by homotopy classes of spheres of weight at most $C_3$
for some constant $C_3>0$.

Let $A$ be a proper submodule of $\pi_2M$. Let $f$ be a solution to
(*) (whose existence is ensured by Lemma~\ref{minimizers}(i)). By the
previous discussion, there exist spheres $f_1,\ldots,f_n$ of weight at
most $C_3$ such that $f$ is in the submodule of $\pi_2M$ generated by
$f_1,\ldots,f_n$. Since $f\not\in A$, at least one of the $f_i$'s, say
$f_1$, is not in $A$. Since $f$ has minimal weight,
$\wt(f)\le\wt(f_1)\le C_3$. Now $f$ is normal by
Lemma~\ref{minimizers}(ii), so by Lemma~\ref{myineq},
$\diam(f)\le C_1$, where $C_1:=C_3^2$.
\end{proof}

We are going to construct inductively a transfinite sequence of
collections $\calS_\lambda$ of maps from $\Ss^2$ to $M$. For some
ordinal $\lambda_0$, the construction will stop, and our system $\calS$
will be obtained by modifying $\calS_{\lambda_0}$.

To start off, set $\calS_0:=\emptyset$ and let $A_0$ be the trivial
$\pi_1M$-submodule of $\pi_2M$. Using Lemma~\ref{diam bound}, we get a
solution $f_1\!:\Ss^2\to M$ to (*) for $A_0$ with diameter at most
$C_1$. We define $\calS_1=\{f_1\}$.

Assuming that $\lambda$ is an ordinal for which $\calS_\lambda$ has been
defined, we let $A_\lambda$ be the $\pi_1M$-submodule of $\pi_2M$
generated by $\calS_\lambda$. As before, Lemma~\ref{diam bound} gives
us a solution $f_{\lambda+1}\!:\Ss^2\to M$ to (*) for $A_\lambda$ with
diameter at most $C_1$ and we put $\calS_{\lambda+1}:=\calS_\lambda
\cup \{f_{\lambda+1}\}$. If $\lambda$ is a limit ordinal, we simply
define $\calS_\lambda$ to be the union of $\calS_\mu$ for all
$\mu<\lambda$.

For some ordinal $\lambda_0$, it occurs that $A_{\lambda_0}=\pi_2M$,
and the construction stops.

\begin{lem}\label{loc finite}
The collection $\calS_{\lambda_0}$ is locally finite.
\end{lem}

\begin{proof}
Let $K$ be a compact subset of $M$. Let $Y$ be a regular neighborhood
of the $(C_1+1)$-neighborhood of $K$. By Kneser-Haken finiteness,
there exists an integer $n>0$ such that in any spherical system in $Y$
of cardinal greater than $n$, one component bounds a ball or two
components are parallel.

Looking for a contradiction, suppose that infinitely many components
of $\calS_{\lambda_0}$ meet $K$. By the diameter bound, infinitely
many components of $\calS_{\lambda_0}$ are contained in $Y$. Let
$f_1,\ldots,f_{n+1}$ be a subcollection of them.

By Lemma~\ref{minimizers}(ii)~and~(iii), they are embeddings or double
covers of projective planes, and pairwise disjoint. Moreover, each of
them is nontrivial in $\pi_2M$ and no two of them are freely homotopic
in $M$. Hence after small homotopies on the double covers of
projective planes, we get a collection of pairwise disjoint
$2$-spheres embedded in $Y$, all homotopically nontrivial and pairwise
nonparallel. This contradiction proves Lemma~\ref{loc finite}.
\end{proof}

Now that we know that $\calS_{\lambda_0}$ is locally finite, we can
perform a perturbation as in the proof of Lemma~\ref{loc finite}
on all components of
$\calS_{\lambda_0}$ at once. This yields a locally finite collection $\calS$ of
pairwise disjoint embedded $2$-spheres that generate $\pi_2M$ as a
$\pi_1M$-module. By Proposition~\ref{top equiv}, $\calS$ splits $M$
into weakly irreducible submanifolds, and the proof of the Main
Theorem is complete.

\begin{concrems}
\item No attempt has been made towards the greatest generality. It
should be straightforward to extend our proof of Theorem~\ref{main} to
nonorientable manifolds, manifolds with boundary, or metrics with mild
singularities such as cone-manifolds. In another direction, the
bounded geometry hypothesis can probably be weakened somewhat.
\item In the course of the proof, we have proven a PL version of the
main theorem. The statement is the same except that ``Riemannian metric
of bounded geometry'' and ``area'' should be replaced by ``triangulation''
and ``weight'' respectively. Note that no bounded geometry hypothesis
is needed in this context.
\item Instead of using PL minimal surfaces, one might want to work
directly with minimal surfaces in the Riemannian manifold $(M,g)$. A
technical problem is that minimizers need not exist because $M$ is
noncompact. This difficulty can be overcome by replacing $g$ by
another metric so that minimizers exist, and using known estimates on
stable minimal surfaces in $3$-manifolds (see e.g.~\cite{schoen}) to
replace Lemma~\ref{myineq}. However, the proof presented here is more
elementary (because existence of PL minimal surfaces is easier to
establish that that of minimal spheres) and we believe the PL version
of the theorem to be interesting in its own right.
\end{concrems}

\def\cprime{$'$}


\begin{thebibliography}{10}

\bibitem{agol:tame}
I.~Agol.
\newblock {Tameness of hyperbolic 3-manifolds}.
\newblock arXiv:math.GT/0405568.

\bibitem{attie:bg}
O.~Attie.
\newblock Quasi-isometry classification of some manifolds of bounded geometry.
\newblock {\em Math. Z.}, 216(4):501--527, 1994.

\bibitem{bt:dehn}
J.~Burillo and J.~Taback.
\newblock Equivalence of geometric and combinatorial {D}ehn functions.
\newblock {\em New York J. Math.}, 8:169--179 (electronic), 2002.

\bibitem{cg:tame}
D.~Calegari and D.~Gabai.
\newblock {Shrinkwrapping and the taming of hyperbolic 3-manifolds}.
\newblock arXiv:math.GT/0407161.

\bibitem{glp:struc}
M.~Gromov, J.~Lafontaine, and P.~Pansu.
\newblock {\em Structures m\'etriques pour les vari\'et\'es riemanniennes}.
\newblock CEDIC, Paris, 1981.

\bibitem{jr:min}
W.~Jaco and J.~H. Rubinstein.
\newblock P{L} minimal surfaces in 3-manifolds.
\newblock {\em J. Differential Geom.}, 27(3):493--524, 1988.

\bibitem{kneser:sphere}
H.~Kneser.
\newblock Geschlossene {F}l\"achen in dreidimensionalen
  {M}annig\-faltig\-keiten.
\newblock {\em Jber. Deutsch. Math.-Verein.}, 38:248--260, 1929.

\bibitem{sm:seifert}
S.~Maillot.
\newblock Open 3-manifolds whose fundamental groups have infinite center, and a
  torus theorem for 3-orbifolds.
\newblock {\em Trans. Amer. Math. Soc.}, 355(11):4595--4638, 2003.

\bibitem{my:embedding}
W.~H. Meeks, III and S.~T. Yau.
\newblock Topology of three-dimensional manifolds and the embedding problems in
  minimal surface theory.
\newblock {\em Ann. of Math. (2)}, 112(3):441--484, 1980.

\bibitem{schoen}
R.~Schoen.
\newblock Estimates for stable minimal surfaces in three-dimensional manifolds.
\newblock In {\em Seminar on minimal submanifolds}, volume 103 of {\em Ann. of
  Math. Stud.}, pages 111--126. Princeton Univ. Press, Princeton, NJ, 1983.

\bibitem{scott:noncompact}
P.~Scott.
\newblock Fundamental groups of non-compact {$3$}-manifolds.
\newblock {\em Proc. London Math. Soc. (3)}, 34(2):303--326, 1977.

\bibitem{st:exotic}
P.~Scott and T.~Tucker.
\newblock Some examples of exotic noncompact $3$-manifolds.
\newblock {\em Quart. J. Math. Oxford Ser. (2)}, 40(160):481--499, 1989.

\bibitem{whitehead:unity}
J.~H.~C. Whitehead.
\newblock A certain open manifold whose group is unity.
\newblock {\em Quart. J. Math. Oxford}, 6:268--279, 1935.

\end{thebibliography}
\end{document}